\documentstyle[12pt]{amsart}
\newtheorem{theorem}{Theorem}[section]

\newtheorem{remark}[theorem]{Remark}

\newtheorem{proposition}[theorem]{Proposition}

\numberwithin{equation}{section}

\begin{document}

\title[The mean curvature flow approach to the symplectic isotopy problem]
{The mean curvature flow approach to the symplectic isotopy
problem}

\author{Xiaoli Han, Jiayu Li}

\thanks{The research was supported by  NSFC }

\address{Academy of Mathematics and Systems Sciences\\
Chinese Academy of Sciences\\ Beijing 100080, P. R. of China.}
\email{hxiaoli@@amss.ac.cn}

\address{Math. Group, The abdus salam ICTP\\ Trieste 34100,
   Italy\\
   and Academy of Mathematics and Systems Sciences\\ Chinese Academy of
Sciences\\ Beijing 100080, P. R. of China. }
\email{lijia@@amss.ac.cn}


\begin{abstract}
Let $M$ be a K\"ahler-Einstein surface with positive scalar
curvature. If the initial surface is sufficiently close to a
holomorphic curve, we show that the mean curvature flow has a
global solution and it converges to a holomorphic
curve.
\end{abstract}

\maketitle

\section{Introduction}

It was conjectured by Tian [T] that every embedded orientable
closed symplectic surface in a compact K\"ahler-Einstein surface
is isotopic to a symplectic minimal surface in a suitable sense.
When K\"ahler-Einstein surfaces are of positive scalar curvature,
this was proved for lower degrees by using pseudo-holomorphic
curves (cf. [ST], [Sh]). Since 1998, we have been trying to study
Tian's conjecture by using the mean curvature flow (c.f. [CT],
[CL1], [CL2]). The mean curvature flow was introduced and
intensively studied by Huisken ([H], [EH]).

The progress has been made during these years.

Suppose that $M$ is a K\"ahler-Einstein surface. Let $\omega$ be
the symplectic form on $M$ and let $J$ be a complex structure
compatible with $\omega$. The Riemannian metric $\langle,\rangle$
on $M$ is defined by
$$
\langle U,V \rangle =\omega(U,JV).
$$
For a compact oriented real surface $\Sigma$ which is smoothly
immersed in $M$, one defines, following [CW], the k\"ahler angle
$\alpha$ of $\Sigma$ in $M$ by
$$\omega|_\Sigma=\cos\alpha d\mu_\Sigma$$ where $d\mu_\Sigma$ is
the area element of $\Sigma$ in the induced metric from
$\langle,\rangle$. We say that $\Sigma$ is a holomorphic curve if
$\cos\alpha \equiv 1$, $\Sigma$ is a Lagrangian surface if
$\cos\alpha \equiv 0$ and $\Sigma$ is a symplectic surface if
$\cos\alpha > 0$.

The evolution equation for the K\"ahler angle was obtained in [CT]
and [CL1] (also see [Wa]). Let $J_{\Sigma_t}$ be an almost complex
structure in a tubular neighborhood of $\Sigma_t$ on $M$ with
\begin{equation}\label{eq2}
\left\{\begin{array}{clcr} J_{\Sigma_t}e_1&=&e_2\\
J_{\Sigma_t}e_2&=&-e_1\\ J_{\Sigma_t}v_1&=&v_2\\
J_{\Sigma_t}v_2&=&-v_1.
\end{array}\right.
\end{equation}
It is not difficult to check (c.f. [CL1]), with $\overline\nabla$
being the covariant derivative of the metric $\langle,\rangle$ on
$M$, that
\begin{equation}\label{djh}
|\overline{\nabla}J_{\Sigma_t}|^2=
|h_{11}^2+h_{12}^1|^2+|h_{21}^2+h_{22}^1|^2
+|h_{12}^2-h_{11}^1|^2+|h_{22}^2-h_{21}^1|^2\geq \frac{1}{2}|H|^2.
\end{equation}
\begin{proposition}([CL1], [CT])
Suppose that $M$ is a K\"ahler-Einstein surface with scalar
curvature $R$. Let $\alpha$ be the K\"ahler angle of the surface
$\Sigma_t$ which evolves by the mean curvature flow. Then
\begin{equation}\label{kef} (\frac{\partial}{\partial t}-\Delta )\cos\alpha =
|\overline{\nabla}J_{\Sigma_t}|^2\cos\alpha +
R\sin^2\alpha\cos\alpha \end{equation} where $R$ is the scalar
curvature of $M$.
\end{proposition}

One corollary of the evolution equation is that, if the initial
surface is symplectic, along the mean curvature flow, at every
time $t$ the surface is still symplectic.

It is proved in [CLT] that, if the scalar curvature of the
K\"ahler-Einstein surface is positive and the mean curvature flow
with initial surface symplectic exists globally, it sub-converges
to a holomorphic curve at infinity outside a finite set of points.

To study the global existence of the mean curvature flow with
initial surface symplectic, we ([CL1], [CL2]) derived a new
monotonicity inequality and proved that the tangent cone at the
first blow-up time is flat which implies that there is no Type I
singularity (also see [Wa]).

In this paper we will show that, {\it if the scalar curvature of
the K\"ahler-Einstein surface is positive and the initial surface
is sufficiently close to a holomorphic curve, the mean curvature
flow has a global solution and it converges to a holomorphic
curve.}

The main new point of this paper is the integral estimate of the
K\"ahler angle
$$\int_{\Sigma_t} \frac{\sin^2\alpha}{\cos\alpha}d\mu_t\leq
\int_{\Sigma_0}\frac{\sin^2\alpha}{\cos\alpha}d\mu_0\exp^{-Rt},$$
which implies the estimation of the mean curvature vector on the
time space,
$$\int_0^T\int_{\Sigma_t}|H|d\mu_t dt\leq
(\int_{\Sigma_0}\frac{\sin^2\alpha}{\cos\alpha}d\mu_0)^{1/2}
\frac{Area(\Sigma_0)^{1/2}}{1-\exp^{-\frac{R}{2}}}.$$

Using the last inequality, we show that the parabolic density
function at any time $t$ along the mean curvature can be dominated
by the one at beginning, then our result follows from White's [W]
regularity theorem.

The result in this paper gives more evidence that the program
should be realized in general.

\section{The main result and its proof}
We first derive an integral estimation of the K\"ahler angle along
the flow, which yields an $L^2$-estimation of the mean curvature
vector on the time space.

\begin{proposition}\label{pro1}
Suppose that $M$ is a K\"ahler-Einstein surface with scalar
curvature $R$. Let $\alpha$ be the K\"ahler angle of the surface
$\Sigma_t$ which evolves by the mean curvature flow. Suppose that
$\cos \alpha(\cdot,0)$ has a positive lower bound. Then
$$\int_{\Sigma_t} \frac{\sin^2\alpha}{\cos\alpha}d\mu_t\leq
C_0\exp^{-Rt}$$ and $$\int^{t+1}_t\int_{\Sigma_t} |H|^2 d\mu_t
dt\leq C_0\exp^{-Rt}$$ where $C_0$ is a constant which depends
only on the initial surface,
$C_0=\int_{\Sigma_0}\frac{\sin^2\alpha}{\cos\alpha}d\mu_0$.
\end{proposition}
{\it Proof:} By applying the parabolic maximum principle to the
evolution equation (\ref{kef}), one concludes that $\cos\alpha$
remains positive as long as the mean curvature flow has a smooth
solution, no matter $R$ is positive, 0 or negative (c.f. [CT],
[CL1], [Wa1]). Moreover, we can obtain that
$$(\frac{\partial}{\partial
t}-\Delta )\frac{1}{\cos\alpha} =-\frac{1}{\cos^2\alpha}
(|\overline{\nabla}J_{\Sigma_t}|^2\cos\alpha +
R\sin^2\alpha\cos\alpha)-\frac{2|\nabla\cos\alpha|^2}{\cos^3\alpha}.$$
Because $\omega$ is closed, we can see that
$$
\int_{\Sigma_t}\cos\alpha d\mu_t = \int_{\Sigma_t}\omega
$$
is constant under the continuous deformation in $t$. In fact, for
any $t_1<t_2$, let $\Omega_{t_1,t_2}=\cup_{t_1\leq t\leq
t_2}\Sigma_t$ be the body spread by $\Sigma_t$ from $t_1$ to
$t_2$, since $\partial \Sigma_{t_1}=\partial
\Sigma_{t_2}=\emptyset$, by Stokes theorem, one has
$$
\int_{\Sigma_{t_2}}\omega -\int_{\Sigma_{t_1}}\omega
=\int_{\Omega_{t_1,t_2}}d\omega =0.
$$

We therefore have
\begin{eqnarray}\label{e1}
\frac{\partial}{\partial
t}\int_{\Sigma_t}\frac{\sin^2\alpha}{\cos\alpha}d\mu_t&=&\frac{\partial}{\partial
t}\int_{\Sigma_t}\frac{1}{\cos\alpha}d\mu_t=\int_{\Sigma_t}(\frac{\partial}{\partial
t}-\Delta)\frac{1}{\cos\alpha}d\mu_t-\int_{\Sigma_t}\frac{|H|^2}{\cos\alpha}d\mu_t\\&=&
-\int_{\Sigma_t}\frac{|\overline\nabla
J|^2}{\cos\alpha}d\mu_t-\int_{\Sigma_t}R\frac{\sin^2\alpha}{\cos\alpha}d\mu_t-\int_{\Sigma}
\frac{|\nabla\cos\alpha|^2}{\cos^3\alpha}d\mu_t-\int_{\Sigma_t}\frac{|H|^2}{\cos\alpha}\nonumber.
\end{eqnarray}
So,
$$\frac{\partial}{\partial t}\int_{\Sigma_t}\frac{\sin^2\alpha}{\cos\alpha}d\mu_t
\leq -R\int_{\Sigma_t}\frac{\sin^2\alpha}{\cos\alpha}d\mu_t.$$ It
is easy to get that
$$\int_{\Sigma_t}\frac{\sin^2\alpha}{\cos\alpha}d\mu_t
\leq\int_{\Sigma_0}\frac{\sin^2\alpha}{\cos\alpha}d\mu_0
\exp^{-Rt}$${\it
i,e,}$$\int_{\Sigma_t}\frac{\sin^2\alpha}{\cos\alpha}d\mu_t \leq
C_0 \exp^{-Rt}.$$

From (\ref{e1}) we know that $$\int_{\Sigma_t}|H|^2 d\mu_t\leq
-\frac{\partial}{\partial
t}\int_{\Sigma_t}\frac{\sin^2\alpha}{\cos\alpha}d\mu_t.$$
Integrating the above inequality from $t$ to $t+1$ we obtain that
\begin{eqnarray*}\int_t^{t+1}\int_{\Sigma_t}|H|^2 d\mu_t
dt&\leq&\int_{\Sigma_t}\frac{\sin^2\alpha}{\cos\alpha} d\mu_t
-\int_{\Sigma_{t+1}}\frac{\sin^2\alpha}{\cos\alpha}
d\mu_t\\&\leq&\int_{\Sigma_t}\frac{\sin^2\alpha}{\cos\alpha}
d\mu_t\leq C_0\exp^{-Rt}.
\end{eqnarray*}
This proves the proposition.
 \hfill Q.E.D.

We derive an $L^1$-estimation of the mean curvature vector on the
time space.
\begin{proposition}\label{pro2}
Suppose that $M$ is a K\"ahler-Einstein surface with positive
scalar curvature $R$. Let $\alpha$ be the K\"ahler angle of the
surface $\Sigma_t$ which evolves by the mean curvature flow.
Suppose that $\cos \alpha(\cdot,0)$ has a positive lower bound.
Then
$$\int_0^T\int_{\Sigma_t}|H|d\mu_t dt\leq (C_0)^{1/2}
\frac{Area(\Sigma_0)^{1/2}}{1-\exp^{-\frac{R}{2}}}.$$
where the constant $C_0$ is defined in Proposition \ref{pro1}.
\end{proposition}
{\it Proof:} We have \begin{eqnarray*}
\int_0^T\int_{\Sigma_t}|H|d\mu_t
dt&=&\sum_{k=0}^{T-1}\int_k^{k+1}\int_{\Sigma_t}|H|d\mu_t
dt\\&\leq&\sum_{k=0}^{T-1}(\int_k^{k+1}\int_{\Sigma_t}|H|^2 d\mu_t
dt)^{1/2}(\int_k^{k+1}area\Sigma_t)^{1/2}\\&\leq&
Area(\Sigma_0)^{1/2}\sum_{k=0}^{T-1}(\int_k^{k+1}\int_{\Sigma_t}|H|^2
d\mu_t dt)^{1/2}\\&\leq&
(C_0)^{1/2}Area(\Sigma_0)^{1/2}\sum_{k=0}^{T-1}\exp^{-\frac{Rk}{2}}\\&\leq&
(C_0)^{1/2}\frac{Area(\Sigma_0)^{1/2}}{1-\exp{\frac{-R}{2}}}.
\end{eqnarray*}
This proves the proposition. \hfill Q.E.D.

Let us recall White's local regularity theorem.

 Let $H (X, X_0 ,t)$ be the backward heat kernel on $R^4$. Define
$$
\rho (X,t)=4\pi (t_0-t)H (X, X_0 ,t)=\frac{1}{4\pi (t_0-t)}\exp
(-\frac{|X-X_0|^2} {4(t_0-t)})
$$
for $t<t_0$. Let $i_M$ be the injective radius of $M^4$. We choose
a cut off function $\phi\in C_0^\infty(B_{2r}(X_0))$ with
$\phi\equiv 1$ in $B_{r}(X_0)$, where $X_0\in M$, $0<2r<i_M$.
Choose a normal coordinates in $B_{2r}(X_0)$ and express $F$ using
the coordinates $(F^1,F^2,F^3,F^4)$ as a surface in $R^4$. The
parabolic density of the mean curvature flow is defined by
$$
\Phi (X_0,t_0,t)=\int_{\Sigma_t}\phi (F )\rho (F ,t)d\mu_t.
$$
The following local regularity theorem was proved by White [W]
(Theorem 3.1 and Theorem 4.1).
\begin{theorem}\label{the1}
There is a positive constant $\varepsilon_0>0$ such that if $$\Phi
(X_0,t_0,t_0-r^2)\leq 1+\varepsilon_0,$$ then the second
fundamental form $A(t)$ of $\Sigma_t$ in $M$ is bounded in $
B_{\frac{r}{2}}(X_0)$, i.e.
$$
\sup_{B_{\frac{r}{2}}\times (t_0-r^2/4,t_0]}|A|\leq C ,
$$
where $C$ is a positive constant depending only on $M$.
\end{theorem}
\begin{remark}\label{rem1}
Since $\Sigma_0$ is smooth, it is well known that
$$\lim_{r\rightarrow 0}\int_{\Sigma_0}\phi(F)\frac{1}{4\pi r^2}
\exp^{\frac{|F-X_0|^2}{4r^2}}d\mu_0=1
$$ for any $X_0\in\Sigma_0$. So we can find sufficiently small
$r_0$ such that
$$\int_{\Sigma_0}\phi(F)\frac{1}{4\pi
r_0^2}\exp^{\frac{|F-X_0|^2}{4r_0^2}}d\mu_0\leq 1+\varepsilon_0/2
$$ for all $X_0\in M$, where $\varepsilon_0$ is constant in White's Theorem.
\end{remark}

Now we state and prove our main result in this paper.
\begin{theorem}\label{the2}
Suppose that $M$ is a K\"ahler-Einstein surface with positive
scalar curvature $R$. Let $\alpha$ be the K\"ahler angle of the
surface $\Sigma_t$ which evolves by the mean curvature flow. Then
there exist sufficiently small constant $\varepsilon_1$ such that
if $C_0\leq\varepsilon_1$ and
$\varepsilon_1/r_0^6\ll\varepsilon_0$, where $C_0$ is defined in
Proposition \ref{pro1}, ${r_0}$ is defined in Remark (\ref{rem1})
and $\varepsilon_0$ is constant in White's theorem, the mean
curvature flow with the initial surface $\Sigma_0$ exists globally
and it converges to a holomorphic curve.
\end{theorem}
{\it Proof:} Fix any positive $T$. By the definition of $\Phi$ we
have
$$\Phi(X_0, t, t-r^2)=\int_{\Sigma_{t-r^2}}\phi(F)
\frac{1}{4\pi r^2}\exp^{-\frac{|F-X_0|^2}{4r^2}}d\mu_{t-r^2}.$$
Differentiating this equation with respect to $t$ we get that
\begin{eqnarray*}
\frac{\partial}{\partial t}\int_{\Sigma_t}\phi(F)\frac{1}{4\pi
r^2}\exp^{-\frac{|F-X_0|^2}{4r^2}}d\mu_t&=&\int_{\Sigma_t}
\bigtriangledown\phi\cdot H \frac{1}{4\pi
r^2}\exp^{-\frac{|F-X_0|^2}{4r^2}}d\mu_t\\&&-\int_{\Sigma_t}\frac{\phi}{8\pi
r^4}\exp^{-\frac{|F-X_0|^2}{4r^2}}\langle F-X_0, H\rangle
d\mu_t\\&&-\int_{\Sigma_t}\frac{\phi}{4\pi
r^2}\exp^{-\frac{|F-X_0|^2}{4r^2}}|H|^2 d\mu_t.
\end{eqnarray*}
Integrating the above equation from $r^2$ to $T$ we get that
\begin{eqnarray*}
\int_{\Sigma_{T-r^2}}\phi(F)\frac{1}{4\pi
r^2}\exp^{-\frac{|F-X_0|^2}{4r^2}}d\mu_t&-&\int_{\Sigma_0}\phi(F)\frac{1}{4\pi
r^2}\exp^{-\frac{|F-X_0|^2}{4r^2}}d\mu_0\\&=&
\int_{r^2}^T\int_{\Sigma_t}\bigtriangledown\phi \cdot H
\frac{1}{4\pi
r^2}\exp^{-\frac{|F-X_0|^2}{4r^2}}d\mu_t\\&&-\int_{r^2}^T\int_{\Sigma_t}\frac{\phi}{8\pi
r^4}\exp^{-\frac{|F-X_0|^2}{4r^2}}\langle F-X_0, H\rangle
d\mu_t\\&&-\int_{r^2}^T\int_{\Sigma_t}\frac{\phi}{4\pi
r^2}\exp^{-\frac{|F-X_0|^2}{4r^2}}|H|^2 d\mu_t.
\end{eqnarray*}

Set $r=r_0$ then we get
\begin{eqnarray*}
\int_{\Sigma_{T-r_0^2}}\phi(F)\frac{1}{4\pi
r_0^2}\exp^{-\frac{|F-X_0|^2}{4r_0^2}}d\mu_t&\leq&\int_{\Sigma_0}\phi(F)\frac{1}{4\pi
r_0^2}\exp^{-\frac{|F-X_0|^2}{4r_0^2}}d\mu_0\\&+&\int_0^T\int_{\Sigma_t}|\bigtriangledown\phi
|| H |\frac{1}{4\pi
r_0^2}\exp^{-\frac{|F-X_0|^2}{4r_0^2}}d\mu_t\\&&-\int_0^T\int_{\Sigma_t}\frac{\phi}{8\pi
r_0^4}\exp^{-\frac{|F-X_0|^2}{4r_0^2}} |F-X_0||H|
d\mu_t\\&&-\int_0^T\int_{\Sigma_t}\frac{\phi}{4\pi
r_0^2}\exp^{-\frac{|F-X_0|^2}{4r_0^2}}|H|^2 d\mu_t.
\end{eqnarray*}
Using Remark (\ref{rem1}), Proposition \ref{pro1}, Proposition
\ref{pro2} and note that $|\bigtriangledown\phi|\leq C$, we obtain
that
\begin{eqnarray*}
\int_{\Sigma_{T-r_0^2}}\phi(F)\frac{1}{4\pi
r_0^2}\exp^{-\frac{|F-X_0|^2}{4r_0^2}}d\mu_t&\leq&1+\varepsilon_0/2+\frac{C}{4\pi
r_0^2}(C_0)^{1/2}\frac{Area(\Sigma_0)^{1/2}}{1-\exp{\frac{-R}{2}}}\\&&+\frac{1}{8\pi
r_0^3}(C_0)^{1/2}\frac{Area(\Sigma_0)^{1/2}}{1-\exp{\frac{-R}{2}}}+\frac{1}{4\pi
r_0^2}C_0\\&\leq&1+\varepsilon_0.
\end{eqnarray*}

Applying White's theorem we obtain an uniform estimate of the
second fundamental form which implies the global existence and
convergence of the mean curvature flow. By Proposition \ref{pro1}
we can see that
$$\int_{\Sigma_t} \frac{\sin^2\alpha}{\cos\alpha}d\mu_t\leq
C_0\exp^{-Rt}.$$ Let $t\rightarrow \infty$ and note that $R>0$, we
get
$$\int_{\Sigma_\infty}
\frac{\sin^2\alpha}{\cos\alpha}d\mu_\infty= 0 .
$$
It follows that $\cos\alpha_\infty=1$, that is $\Sigma_\infty$ is
a holomorphic curve. This proves the theorem. \hfill Q.E.D.

 \vspace{.2in}
\begin{center}
{\large\bf REFERENCES}
\end{center}
\footnotesize
\begin{description}
\item [{[CL1]}] {J.Chen and J.Li, Mean curvature flow of surface
in 4-manifolds, Adv. Math. {\bf 163} (2001), 287-309.}

\item[{[CL2]}] {J. Chen and J. Li, Singularities of codimension
two mean curvature flow of symplectic surfaces, preprint .}

\item[{[CLT]}] {J. Chen, J. Li and G. Tian, Two-dimensional graphs
moving by mean curvature flow, Acta Math. Sinica, English Series
Vol.{\bf 18}, No.2 (2002), 209-224.}

\item[{[CT]}] {J. Chen and G. Tian, Moving symplectic curves in
K\"ahler-Einstein surfaces, Acta Math. Sinica, English Series,
{\bf 16} (4), (2000), 541-548.}

\item[{[CW]}] {S.S.Chen and R.Wolfson, Minimal surfaces by moving
frames, Amer. J. Math. {\bf 105} (1983), 59-83.}

\item[{[EH]}] {K. Ecker and G. Huisken, {\it Mean curvature
evolution of entire graphs}, Ann. Math. {\bf 130} (1989),
453-471.}

\item[{[H]}] {G. Huisken, Contracting convex hypersurfaces in
Riemannian manifolds by their mean curvature, Invent. math. {\bf
84} (1986), 463-480.}

\item[{[Sh]}] {V. Shevchishin, Pseudo-holomorphic curves and the
symplectic isotopy problem, preprint.}

\item[{[ST]}] {B. Siebert and G. Tian, Holomorphy of genus two
 Lefschetz fibration, preprint.}

\item[{[T]}] {G. Tian, Symplectic isotopy in four dimension, First
International Congress of Chinese Mathematicians (Beijing, 1998),
143-147, AMS/IP Stud. Adv. Math., 20, Amer. Math. Soc.,
Providence, RI, 2001. }

\item[{[W]}] {B.White, A local regularity theorem for classical
mean curvature flow, Preprint, Stanford.}

\item[{[Wa]}] {M.T.Wang, Mean curvature flow of surfaces in
Einstein four manifolds, J. Diff. Geom. {\bf 57} (2001), 301-338.}
\end{description}
\end{document}